\documentclass{amsart}
\usepackage{latexsym,amssymb,graphicx}
\title[Simplicity of the Reduced $C^{\ast}$-Algebras of Coxeter
Groups]{Simplicity of the Reduced $C^{\ast}$-Algebras of Certain Coxeter Groups}
\author[Gero Fendler]{Gero Fendler \\Naturwissenschaftlich--Technische
  Fakult\"at I\\Fachrichtung 6.1 Mathematik \\
Universit\"at des Saarlandes \\ 
}
\date{}
\newcommand{\NN}{\ensuremath{\mathbb{N}}}
\newcommand{\ZZ}{\ensuremath{\mathbb{Z}}}
\newcommand{\CC}{\ensuremath{\mathbb{C}}}
\newcommand{\RR}{\ensuremath{\mathbb{R}}}
\newcommand{\HH}{\ensuremath{\mathbb{H}}}
\newcommand{\ii}{{\rm i}}
\newcommand{\norm}[2][]{\| \, {#2} \,\|_{#1}}

\newcommand{\card}[1][\#]{\protect{#1}}
\newcommand{\Gl}{{\rm Gl}}
\newcommand{\sign}{\ensuremath{\mbox{\rm sign}}}
\newcommand{\la}[1]{{\rm \bf l}(#1)}
\newcommand{\e}{\mathbf{e}}
\newcommand{\Sym}{\ensuremath{\mbox{\rm Sym}}}
\newcommand{\Aut}{{\rm Aut}}
\newcommand{\spec}[1][.]{{\rm Sp}(#1)}
\newcommand{\id}{\ensuremath{\mbox{\rm Id}}}
\newtheorem{thm}{Theorem}
\newtheorem{cor}{Corollary}
\newtheorem{lem}{Lemma}
\newtheorem{rem}{Remark}
\newtheorem{prop}{Proposition}
\begin{document}
\begin{abstract}
{Let $(G,S)$ be a finitely generated Coxeter group,
such that the Coxeter system is indecomposable and the
canonical bilinear form is indefinite but non-degenerate.
We show that the reduced $C^{\ast}$-algebra of $G$
is simple with unique normalised
trace.
\par
For an arbitrary finitely generated Coxeter group we prove the 
validity of a Haagerup inequality:
There exist constants $C>0$ and $\Lambda\in \NN$
such that for a function $f\in l^2(G)$ supported on elements
of length $n$ with respect to the generating set $S$:
\[\norm[2]{f\ast h} \leq C(n+1)^{\frac{3}{2}{\Lambda}} \norm[2]{f}
\norm[2]{h},\quad \forall h\in l^2(G).\]
} 
\end{abstract}
\maketitle
\vfill
\section{Introduction}
For a discrete group $G$ we denote $l^2(G)$ the Hilbert space 
of all square summable complex functions on $G$ and 
$B(l^2(G))$ the von~Neumann algebra of all bounded operators on $l^2(G)$.
\par
The group $G$ acts on $l^2(G)$ by the left regular representation:
\begin{eqnarray*}
\lambda(g) f(h) &=& f(g^{-1}h),\quad g,h \in G,\; f\in l^2(G).
\end{eqnarray*}
The reduced (or (left) regular) $C^{\ast}$-algebra $C^{\ast}_r (G) $  of 
$G$ is the operator norm closure
of the linear span of the set of operators $\{\lambda(g):\; g\in G\}$.
We often think of its elements as certain $l^2$ functions on $G$.
The natural normalised trace on this algebra is then just evaluation of 
a function at the group identity.\par
When $G$ is a non-abelian free group on two generators 
Powers~\cite{Powers_75} showed that $C^{\ast}_r (G)$ is a simple $C^{\ast}$
algebra, that is, it contains no non-trivial twosided ideals.
Since then this result has been generalised to various groups and  extended
to (reduced) cross products
by several authors (see \cite{Akemann_81}, \cite{Bedos_91},
 \cite{BeCoHa95},  \cite{BocaNitica_88}, \cite{Harpe_85}, 
\cite{HarpeSkandalis}, \cite{PaSal_79}).\par
On the other hand,
when $G$ is an amenable group, or contains amenable normal subgroups,
then the kernel of the trivial representation, respectively the
kernel of the representation induced from the trivial representation
of an amenable  normal subgroup, will be a non-trivial twosided ideal in
$C^{\ast}_r (G)$. These facts makes it reasonable
that the question of simplicity of $C^{\ast}_r (G)$ is
related to the Tits alternative for linear groups (a linear group
either is amenable or contains non-abelian free subgroups).
\par
We shall in this note consider finitely generated Coxeter groups
for which de~la~Harpe~\cite{Harpe_87} gave an elaboration on the Tits alternative
and we shall show that a finitely generated infinite Coxeter group
either contains a normal solvable (even nilpotent) subgroup
or has a simple reduced $C^{\ast}$-algebra.
\par
For the geometric representation $\sigma$ of a Coxeter group $(G,S)$,
with $\card{S}<\infty$, on $E =\mathbb{R}^{S}$ 
we adopt the usual notation of~\cite{Bourbaki_68}.
Assume that $(G,S)$ is indecomposable. 
The canonical $\sigma(G)$-invariant bilinear form $B$ can be
strictly positive definite, positive
semidefinite, or  strictly indefinite.
In these cases, respectively, the group is finite, an affine Coxeter group and
hence amenable, or non-amenable~\cite{Harpe_87}.
% in which case $G$ is finite, positive
%semidefinite, $G$ is an affine Coxeter group hence amenable,
%or strictly indefinite, $G$ being non-amenable then~\cite{Harpe_87}.
In this last case it might occur, that $B$ is degenerate. The
orthogonal $E^0$ of $E$ for $B$ then is pointwise fixed by all 
$\sigma(g)$, $ g\in G$,  and the kernel of the representation
$\hat{\sigma} : G \to \Gl(\hat{E})$ induced on the quotient
$ \hat{E}=E/E^0$ is a non-trivial nilpotent normal subgroup.
(Choosing an appropriate basis for $E$ it is easily seen to be mapped
by $\sigma$
into a group of unipotent matrices.)
\par
If $(G,S)$ is a decomposable Coxeter system, then $G$ can be written as a
direct product of indecomposable Coxeter groups. The reduced
$C^{\ast}$-algebra of $G$
is the spatial tensor product of the reduced $C^{\ast}$-algebras of the
factors. This spatial tensor product is known to be simple if and only if each factor is a
simple $C^{\ast}$-algebra~\cite{Tak_64}.
\par
Hence, we
shall always assume that $B$ is indecomposable and strictly indefinite but
non-degenerate. In the development of the arguments we shall see that under
these conditions a Coxeter group is an icc-group (conjugacy classes of
elements different from the identity are infinite) and that the normalised
trace on $C^{\ast}_r (G)$ is unique. (For a decomposable Coxeter group an
argument  similar to the one in the last
paragraph shows that the trace is unique if and only if the trace is unique 
for each indecomposable factor, see~\cite{BekkaHarpe_2000}.)\\
\hfill
\newlength{\foo}
\settowidth{\foo}{\setlength{\unitlength}{2947sp}
\begin{picture}(2315,1681)(1748,-1464)
\end{picture}
}%
\newlength{\fOO}
\setlength{\fOO}{\textwidth}
\addtolength{\fOO}{-\foo}
\setlength{\foo}{\parindent}
\parbox[b]{\fOO}{%
\hspace*{\foo}One might think that Coxeter groups of the above kind are
Gromov hyperbolic and arguments like in \cite{Harpe_88} combined with 
\cite{BekkaHarpe_2000} would allow to prove simplicity of the reduced 
$C^{\ast}$-algebra, but the group with the Coxeter graph in the figure 
does not contain a finite index Gromov hyperbolic subgroup.\hfill}
\hfill
\setlength{\unitlength}{2447sp}%
\begingroup\makeatletter\ifx\SetFigFont\undefined%
\gdef\SetFigFont#1#2#3#4#5{%
  \reset@font\fontsize{#1}{#2pt}%
  \fontfamily{#3}\fontseries{#4}\fontshape{#5}%
  \selectfont}%
\fi\endgroup%
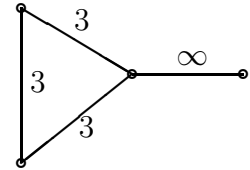
\begin{picture}(2315,1681)(1748,-1500)
\thicklines
\put(1801,164){\circle{90}}
\put(2926,-511){\circle{90}}
\put(4051,-511){\circle{90}}
\put(1801,-1411){\circle{90}}
\put(1801,-1411){\line( 0, 1){1575}}
\put(1801,164){\line( 5,-3){1125}}
\put(2926,-511){\line(-5,-4){1125}}
\put(1801,-1411){\line( 0, 1){  0}}
\put(2926,-511){\line( 1, 0){1125}}
\put(2341,-61){\makebox(0,0)[lb]{\smash{\SetFigFont{12}{14.4}{\rmdefault}{\mddefault}{\updefault}3}}}
\put(1891,-691){\makebox(0,0)[lb]{\smash{\SetFigFont{12}{14.4}{\rmdefault}{\mddefault}{\updefault}3}}}
\put(2386,-1141){\makebox(0,0)[lb]{\smash{\SetFigFont{12}{14.4}{\rmdefault}{\mddefault}{\updefault}3}}}
\put(3376,-421){\makebox(0,0)[lb]{\smash{\SetFigFont{12}{14.4}{\rmdefault}{\mddefault}{\updefault}$\infty$}}}
\end{picture}
\par
Gromov hyperbolic Coxeter groups have been characterized by
Moussong in \cite{Moussong}, as those Coxeter groups $(G,S)$, 
which do not contain two infinite commuting parabolic subgroups, and 
further have the property that
no subset $T\subset S$ generates a parabolic subgroup $(G_T,T)$, which is an
affine Coxeter group of rank at least $3$.
\par
When $B$ has signature $(n-1,1)$ it can happen that $G$ is a
hyperbolic Coxeter group (in the classical sense). Then it is a lattice in 
the real Lie group
$O(n-1,1)$, hence Zariski dense in it. A theorem of Bekka, Cowling and de~la~Harpe
\cite{BeCoHa95} then applies. We do not know whether 
a Coxeter group  of the kind considered here is always Zariski dense
in some simple real Lie group\footnote{This has recently been 
esablished by Benoist and de la Harpe in: Adh{«e}rence de Zariski de 
groups de Coxeter, preprint}. 
\par
To prove simplicity of the reduced $C^{\ast}$-algebra we have to deal
with the combinatorics in $G$. As a byproduct we obtain a
Haagerup inequality, valid for all finitely generated Coxeter groups:
\par
There exist constants $C>0$ and $\Lambda\in \NN$
such that for a function $f$ supported on elements
of word-length $n$ with respect to the generating set $S$:
\[\norm{\lambda(f)} \leq C(n+1)^{\frac{3}{2}{\Lambda}} \norm[2]{f}.\]
\par
The constant $\Lambda$ in this inequality can be obtained
in terms of the geometrical representation of $(G,S)$.
Examples show that it is not best possible. We conjecture
that the optimal constant is just the virtual cohomological
dimension of $G$ and refer the reader to \cite{bestvina93} for
motivation. 
\par
We thank the referee for his comments, which improved the presentation.
\section{Trees}
In this section we shall define certain trees on which a
finite index torsion free normal subgroup $\Gamma$ of the Coxeter group
acts by simplicial automorphisms of the trees. As far as trees are concerned
we use the standard notation of~\cite{Se}, for the existence of a  
finite index torsion free normal subgroup see chap.V, \S 4, ex.9 
of~\cite{Bourbaki_68}. 
We shall show further
that the action of $\Gamma$ on the product of those trees is free.
Our construction is similar to that of
Januszkiewicz~\cite{Januszkiewicz_99}. For the readers convenience we
shall work with the classical Tits cone 
$U$ and the transposed geometrical action $\sigma^{\ast}$ of $G$ on it.
Let us introduce a little notation and recall some facts.
\par
The word-length of $g\in G$, with respect to the generating set $S$, is
defined as
$\la{g}=\inf \{ n : g=s_1\cdot  \ldots \cdot  s_n ,\; s_1,\ldots,s_n \in S\}$.
We denote by $T=\{g^{-1}sg : g\in G, s\in S\}$ the set of reflections of $G$.
Let for $g\in G$ $N_g=\{t\in T: \la{tg}<\la{g}\}$.
With these notations, for $g,h\in G$:
$$ \la{g}=\card\{t\in T: \la{tg}< \la{g}\}.$$
Moreover, see~\cite{BozJanSpatz_88},
$$\la{g^{-1}h}=\card N_g \triangle N_h.$$
\par
We decompose the set of reflections $T \subset G$
in disjoint $\Gamma$-orbits with respect to conjugation:
\begin{eqnarray}%
\label{eq:orbits}
T&=&T_1 \dot{\cup} T_2 \dot{\cup} \ldots \dot{\cup} T_{\Lambda}.
\end{eqnarray}
To $t\in T$ denote by $M_t$ the hyperplane in $E^{\ast}$
fixed by $\sigma^{\ast}(t)$ and call it the mirror of $t$.
For $i \in \{1, \ldots , \Lambda\}$ we define a graph $\mathcal{T}_i$
as follows:
The vertices are the connected components of 
$U \setminus (\bigcup_{t\in T_i}M_t)$
and two such vertices are connected by an edge if, as connected components,
they are separated by just one mirror.
\begin{lem}
\label{lem:tree}
The above defined graph is a tree.
\end{lem}
\begin{proof}
We have to show that a closed path in $\mathcal{T}_i$ contains backtracking.\par
Let $C_0, C_1,\ldots,C_n = C_0, n \geq 1$ be the sequence
of vertices of a non-trivial closed path.
Choose points $c_i \in C_i$, where we may assume $c_0=c_n$,
 and elements $e_1,\ldots , e_n$
of $E$, considered as 
functionals on $E^{\ast}$, with 
$e_i(c_0) < 0$ for all $i\in \{1,\ldots,n\}$,
in such a way that 
$e_i$ vanishes on the hyperplane, which defines the edge $\{C_{i-1},C_i\}$.
We may and do assume that $e_i=e_j$ if the defining edges $\{C_{i-1},C_i\}$
and $\{C_{j-1},C_j\}$ are equal.
\par
Consider the function $f:\{0,\ldots,n\} \mapsto \ZZ$ defined by
\[
f(i) = \sum_{j=1}^{n} \sign \ e_j(c_i).
\]
It fulfils
$f(0) = f(n) =  -n$, and $f(1)  = f(n-1)  =  -n +2$.
\par
Since $f(i+1) \in \{f(i)+2, f(i)-2\}$, we find
$i_0$ such that $f(i_0) > f(i_0+1)=f(i_0-1)$.
Hence $e_{i_0}(c_{i_0}) >0, e_{i_0+1}(c_{i_0}) >0$,
$e_{i_0}(c_{i_{0-1}}) <0$, $e_{i_0+1}(c_{i_0+1}) <0$,
and of course $e_{i_0}(c_0) <0$, $e_{i_0+1}(c_0)<0$.
Since $U$ is convex the hyperplanes $e_{i_0}=0$ and
$e_{i_0+1}=0$ must intersect inside $U$ and we conclude
from Lemma 3 of \cite{fendler99a} that they coincide.
We have found a backtracking.
\end{proof}
\par
The Coxeter group $G$ acts via $\sigma^{\ast}$
on the chamber system $\mathcal{C}$ defined from the mirrors on $U$.
Our basic reference here is  \cite{Bourbaki_68},
but one should also compare with \cite{Ronan89}.
Moreover, two points $x,y$ are separated by a mirror $M_t$ if
and only if, for $g\in G$, $\sigma^{\ast}(g)x$ and $\sigma^{\ast}(g)y$
are separated by $M_{gtg^{-1}}$. Since we defined the trees $\mathcal{T}_i$
with respect to a $\Gamma$-orbit in $T$ we have:
\begin{lem}%
\label{lem:treeauto}
The contragradient representation  $\sigma^{\ast}$ induces an action 
of  $\Gamma$ on $\mathcal{T}_i$ by automorphisms of the tree.
\end{lem}
\begin{proof}
First we show that the action of $\Gamma$ is well defined.
Let $\gamma \in \Gamma$ and a component $C$ be given.
For $c_0,c_1 \in C$ we have to show that $\sigma^{\ast}(\gamma)c_0$
and $\sigma^{\ast}(\gamma)c_1$ are not separated by a mirror of 
a reflection in $T_i$. Indeed, assume that $\sigma^{\ast}(\gamma)c_0$
and $\sigma^{\ast}(\gamma)c_1$ are separated
by $M_t$, then by the remark
before the lemma $M_{\gamma^{-1} t \gamma}$ would separate $c_0$ 
and $c_1$. 
\par
Let $\gamma \in \Gamma $ be given
If $C_0$ and $C_1$ are connected by an edge,
then there exist exactly one $t\in T_i$ such that
$\overline{C_0} \cap M_t \neq \emptyset$ and $\overline{C_1} \cap M_t \neq %
\emptyset$. Clearly this is the case if and only if 
$\overline{\sigma^{\ast}(\gamma)C_0} \cap M_{\gamma{t}\gamma^{-1}} \neq \emptyset$ and 
$\overline{\sigma^{\ast}(\gamma)C_1} \cap M_{\gamma{t}\gamma^{-1}} \neq %
\emptyset $. Since $\gamma{t}\gamma^{-1}\in T_i$ exactly if $t\in T_i$ we 
are done.
\end{proof}
\par 
We consider the product
\begin{eqnarray}%
\label{eq:prod}
\mathcal{G}&=&\mathcal{T}_1 \times \ldots \times \mathcal{T}_{\Lambda}
\end{eqnarray}
as a product of chamber systems (see  \cite{Ronan89} p. 2). 
On the vertices $V_\mathcal{G}$ of $\mathcal{G}$ 
we use the metric:
\[d^1(x,y) = \sum_{i=1}^{\Lambda} d_i(x_i,y_i),\]
where $x=(x_1,\ldots, x_{\Lambda}),y=(y_1,\ldots,
y_{\Lambda})\in V_\mathcal{G}$. The action of $\Gamma$ on $V_\mathcal{G}$
is isometric with respect to this metric.
\begin{lem}
$\Gamma$ acts freely on the vertices of $\mathcal{G}$ without bounded orbit.
Moreover, no non-trivial subgroup of $\Gamma$ has a bounded orbit.
\end{lem}
\begin{proof}
Denoting $C_0$ the fundamental chamber in $U$ we have
the injection $g \mapsto \sigma^{\ast}(g)C_0$ of $G$ 
onto the chambers of  $\mathcal{C}$.
If $[C]_i$ denotes the the connected component of 
the chamber $C$ in $U \setminus (\bigcup_{t\in T_i}M_t)$,
we obtain a map of the chambers of
$\mathcal{C}$ into the set of vertices of the product $\mathcal{G}$:
$[.]:C \mapsto ([C]_1, \ldots , [C]_{\Lambda})$.
It is an injection, since two chambers $C,C'$ in $\mathcal{C}$
are different, if they are separated by a mirror, say by $M_t$ where $t\in T$.
If so then there is $j_0 \in \{1,\ldots, \Lambda\}$ with $t\in T_{j_0}$.
Whence $[C]_{j_0} \neq [C']_{j_0}$. 
\par
The composition of these two maps defines an embedding of
$G$ into the vertices
of $\mathcal{G}$ and the action of $\Gamma$ on this subset is free, 
since it is just the transfered left multiplication in the group $G$.
Moreover, no non-trivial subgroup of $\Gamma$ has a bounded orbit.
\par
To see this we note first, that the injection
$g \mapsto [\sigma^{\ast}(g)C_0]$
is an isometry from $G$ endowed with
the left invariant distance coming from the word-length
with respect to the generating set $S$ into the vertices of $\mathcal{G}$
endowed with the metric $d^1$.
This latter holds true because for $g,h \in G$:
\begin{eqnarray*}
\la{g^{-1}h}&=& \card{N_g \triangle N_h} \\
&=& \card \{ M_t : M_t \mbox{ separates } \sigma^{\ast}(g)C_0 \mbox{ from }% 
\sigma^{\ast}(h)C_0 \}\\
&=& d^1([\sigma^{\ast}(g)C_0],[\sigma^{\ast}(h)C_0]).
\end{eqnarray*}
\par
So, if $ [\sigma^{\ast}(\gamma^n)C_0], n\in \ZZ$ were bounded in
$\mathcal{G}$ then the set of mirrors which separate
$C_0$ from a chamber in $\bigcup_{n\in \ZZ}\sigma^{\ast}(\gamma^n)C_0$
would have finite cardinality. Hence $\sup_{n\in \ZZ} \la{\gamma^n} < \infty$.
We infer that the set $\{\gamma^n :n\in \ZZ \}$ would be finite.
This is a contradiction to the fact that $\Gamma$ is torsion free.
\par
Now, if $x\in \mathcal{G}$ has a stabiliser $\Gamma_x \subset \Gamma$
then a vertex $w=[\sigma^{\ast}(g)C_0]$ in the image of 
$G$ in $ \mathcal{G}$ would have a
bounded $\Gamma_x$-orbit, since $\Gamma$ acts by isometries.
This follows from
\[ d^1(\gamma w, w) \leq d^1(\gamma w, \gamma x) + d^1(\gamma x, x) + d^1(x,w)
\leq 2 d^1(x,w), \quad \forall\gamma \in \Gamma_x. \] 
Hence $\Gamma_x=\{ \e \}$.
\end{proof}
\section{The Action on the Trees}
In this section we shall collect some
auxiliary results for later use.
\begin{lem}%
\label{lem:trans}%
If $t_1,t_2 \in T$
are reflections such that the corresponding edges are distinct 
but in the same tree
then $t_1t_2$ acts as a 
translation on this tree.
\end{lem}
\begin{proof}
First note that $t_1$ and $t_2$ are $\Gamma$ conjugate, $\gamma^{-1}t_1\gamma=
t_2$ say, since their edges belong to the
same tree, $\mathcal{T}_i$ say.
Then, $t_1t_2= t_1\gamma^{-1}t_1\gamma \in \Gamma$.
\par 
An oriented line segment in the Tits-cone,
from a point $v\in M_{t_1}$ to its image
$\sigma^{\ast}(t_2)v$ is just reversed by $\sigma^{\ast}(t_2)$.
Since the edges are distinct, this implies that this segment is non-trivial.
Since $\sigma^{\ast}(t_1)$ maps this line segment to one 
adjacent (both segments contain $v$), but differently oriented, 
we conclude that 
the composition  $\sigma^{\ast}(t_1)\sigma^{\ast}(t_2)$
maps the original segment to a coherently oriented one.
\par
Its image, under $\sigma^{\ast}(t_1t_2)$,
and  the line segment itself
can be connected to a coherently oriented broken line in the cone.
The mirrors crossed by the line segment and those
crossed by its $\sigma^{\ast}(t_1t_2)$-image are separated
by $M_{t_1}$.
Hence, in
$\mathcal{T}_i$, this broken line defines a coherently oriented geodesic.
\end{proof}
\par
The edges of one of the trees $\mathcal{T}_i$ (identified with the set of
reflections $T_i$), as a $\Gamma$-orbit of a 
reflection, generate a subgroup in $G$. 
By a theorem independently proved by  Deodhar \cite{Deod89}
and Dyer \cite{dyer90}  this subgroup is
itself a Coxeter group. Clearly this subgroup is
normalised by $\Gamma$, but in general we can not expect that
all its reflections
are contained in $T_i$.
\par
This improves for subgroups generated by a $G$-conjugation invariant
set of reflections:
\begin{lem}%
\label{lem:reflectionsubgroup}
Let $T'\subset T$ be a set of reflections of $G$, invariant under
conjugation. Let $W'$ denote the  subgroup generated by $T'$ in $G$.
It is,
with respect to a subset $S'\subset T'$, a Coxeter group, normal in $G$,
and
its set of reflections coincides with $T'$.
\end{lem}
\begin{proof}
From the theorem of \cite{Deod89}, or rather from Step~1 of its proof,
(see also Theorem~3.4 and Corollary ~3.11 in \cite{dyer90})
it is clear that $(W',S')$ is a Coxeter system for some set $S'\subset T'$.
A reflection in $W'$ is conjugate, by an element of $W'$, to 
a reflection in $S'$. Since $T'$ is $G$-conjugation invariant,
any reflection of $W'$ is in $T'$.
The other assertions of the lemma are immediate.
\end{proof}
\par
Now we shall view the set of edges of the product of trees (\ref{eq:prod}) 
as a fiber bundle
$p : \mbox{edges}(\mathcal{T}_1 \times \ldots \times \mathcal{T}_{\Lambda})\rightarrow%
\{1, \ldots, \Lambda\}$
with base space $\{1, \ldots, \Lambda\}$. Indeed, two vertices $x=(x_1,\ldots,
x_{\Lambda})$, $y=(y_1,\ldots,y_{\Lambda})$ are connected by an edge,
call it  $e(x,y)$, if for one $j\in \{1,\ldots, \Lambda\}$ the vertices $x_j$
and $y_j$ are connected by an edge in $\mathcal{T}_j$ and for all $i\neq j$
we have
$x_i=y_i$. We define $p(e(x,y))=j$.
\par
Since $\Gamma$ leaves the fibers invariant we obtain
an action of $G/\Gamma$ by permutations of $\{1, \ldots, \Lambda\}$, which we
denote by
$\pi :  G/\Gamma \rightarrow  \Sym_{\Lambda}$.
If $\mathcal{O} \subset \{1, \ldots, \Lambda\}$ is a $\pi(G/\Gamma)$-orbit,
then, by Lemma~\ref{lem:reflectionsubgroup},
the edges of $p^{-1}(\mathcal{O})$ are the reflections of
a Coxeter group $W_\mathcal{O}\lhd G$.
\begin{lem}%
\label{lem:morph}
For $\ i \in \{1, \ldots, \Lambda\}$ and $g \in G$ 
$\sigma^{\ast}(g)$ induces a morphism of trees:
\begin{eqnarray*}
g : \mathcal{T}_i & \rightarrow & \mathcal{T}_{\pi(\dot{g})(i)}.
\end{eqnarray*}
Here $g \mapsto \dot{g}$ denotes the quotient morphism 
$G \rightarrow G/\Gamma$.
\end{lem}
\begin{proof}
An edge of $\mathcal{T}_i$ is the mirror $M_t$ of some
reflection $t\in T_i$. Its image under $\sigma^{\ast}(g)$
is the mirror of $gtg^{-1}\in T_{\pi(\dot{g})(i)}$. Hence it defines an edge
in $ \mathcal{T}_{\pi(\dot{g})(i)}$.
\par
Given a component $C$ of $U\setminus \{M_t: t\in T_i\}$ some
$c_0,c_1\in C$ would have images $\sigma^{\ast}(g)c_0,\sigma^{\ast}(g)c_1$
in different components of $U\setminus \{M_t: t\in T_{\pi(\dot{g})(i)}\}$
if there is a mirror of some $t'\in T_{\pi(\dot{g})(i)}$ separating the images.
But then $g^{-1}t'g \in T_i$ would have a mirror separating $c_0$ and $c_1$.
This contradiction shows that $\sigma^{\ast}(g)$
defines a map from the vertices of $\mathcal{T}_i$ to those of
$\mathcal{T}_{\pi(\dot{g})(i)}$.
\par
If $C_0$ and $C_1$ are connected by an edge in $\mathcal{T}_i$,
then there exists exactly one $t\in T_i$ such that
$\overline{C_0} \cap M_t \neq \emptyset$ and $\overline{C_1} \cap M_t \neq %
\emptyset$. Clearly this is the case if and only if 
$\overline{\sigma^{\ast}(g)C_0} \cap M_{g{t}g^{-1}} \neq \emptyset$ and 
$\overline{\sigma^{\ast}(g)C_1} \cap M_{g{t}g^{-1}} \neq %
\emptyset $. Since $g{t}g^{-1}\in T_{\pi(\dot{g})(i)}$ 
exactly if $t\in T_i$ we 
are done.
\end{proof}
\section{A Haagerup Inequality}
First, following Rammage, Robertson and Steger
\cite{RaRoSt98} we shall prove a Haager\-up inequality
for the torsion free subgroup $\Gamma$ of $G$. Then we shall
apply a theorem of Jolissaint \cite{Jolissaint90} to the group extension
$0\rightarrow \Gamma  \rightarrow G \rightarrow G/{\Gamma} \rightarrow 0$.
\par
We consider the product of trees $\mathcal{G}$ as a building of type
$\tilde{A}_1 \times \ldots \times \tilde{A}_{1}$. Its apartments 
are $\Lambda$-dimensional euclidian spaces tessellated by unit cubes.
We have a shape defined on pairs of vertices
\[
\sigma : V_{\mathcal{G}}\times V_{\mathcal{G}} \rightarrow \ZZ_+\times \ldots%
\times \ZZ_+
\]
by
\[
\sigma(u,w)=(d_1(u_1,w_1),\ldots ,d_{\Lambda}(u_{\Lambda},w_{\Lambda})).
\]
It is clear from Lemma \ref{lem:treeauto} that the action of $\Gamma$ 
is shape preserving and we define a shape on $\Gamma$ by fixing a vertex
$v_0\in V_{\mathcal{G}}$: $\sigma(\gamma)=\sigma(v_0,\gamma v_0)$.
Let $p(n_1,\ldots,n_{\Lambda})= \prod_{i=1}^\Lambda(n_i+1)$. The following
theorem can be proved almost verbatim as the $\tilde{A}_1\times\tilde{A}_1$ case
of Theorem 1.1 of \cite{RaRoSt98}:
\begin{thm}
If $h\in l^2(\Gamma)$ is supported
on elements of shape $(n_1, \ldots , n_{\Lambda})$  then for 
$f \in l^2(\Gamma)$:
\[ \norm[2]{f\ast h} \leq p(n_1,\ldots,n_{\Lambda})\norm[2]{f}\norm[2]{h}.\]
\end{thm}
\begin{cor}
Let $(G,S)$ be a Coxeter group. There exist constants $C>0$ and $\Lambda\in \NN$
such that for a function $h\in l^2(G)$ supported on elements
of length $n$, for all $f\in l^2(G)$:
\[\norm[2]{f\ast h} \leq C(n+1)^{\frac{3}{2}{\Lambda}} \norm[2]{f} \norm[2]{h}.\]
\end{cor}
\begin{proof}
Let $\Gamma$ be a torsion free, normal  subgroup of finite index in
$G$ and denote by $\Lambda$ the
cardinality of distinct conjugation orbits of $\Gamma$ on the set
of reflections of $G$. Then, since the length of an element of $\Gamma$
is just the sum of the components of its shape, we obtain that the set
of elements of length $n$ decomposes in less than 
$k=(n+1)^{\Lambda}$ sets of elements of different shapes.
Obviously $p(\sigma(\gamma))\leq (\la{\gamma}+1)^{\Lambda}$. Hence for 
$h \in l^2(\Gamma)$ with support in elements of length $n$:
\begin{eqnarray*}
\norm[2]{f\ast h} &=& \norm[2]{\sum_{j=1}^k f\ast h_j} \\
&\leq&
(n+1)^{\Lambda}\sum_{j=1}^k\norm[2]{f}\norm[2]{h_j}\\
& \leq&%
(n+1)^{\Lambda}\sqrt{k}\norm[2]{f}\norm[2]{h},
\end{eqnarray*}
where  $h = \sum_{j=1}^k h_j$ is the orthogonal decomposition of $h$ 
into functions $h_j$ supported on elements of the same shape.
\par
Since $\Gamma$ is of finite index in $G$, we may apply
Lemma 2.1.2 of \cite{Jolissaint90}.
\end{proof}
\section{Free Subgroups}
As before let $\Gamma$ be a torsion free subgroup of
finite index in the Coxeter group $G$ and $\mathcal{T}_1, \ldots, \mathcal{T}_\Lambda$
the associated trees.
\par
It is obvious from the definition of the trees, that for each
of them the action
of $\Gamma$ on the set of its edges is transitive.
Hence $\Gamma \backslash \mathcal{T}_i$ is either a simple loop
or a single edge with two endpoints, depending on whether
$\Gamma$ has one or two orbits on the set of vertices.
\par
\begin{lem}%
\label{lem:translation}
For $\gamma \in \Gamma$ there exists a tree,
among $\mathcal{T}_1, \ldots, \mathcal{T}_\Lambda$,
on which
$\gamma$ acts as a translation.
\end{lem}
\begin{proof}
Denote in $\mathcal{T}_i$ by $v_i(\e)$ the vertex defined by the
equivalence class of the group identity.
Since $\la{\gamma^n}\rightarrow \infty$, as $n\rightarrow \infty$,
the formula
\begin{eqnarray*}
\la{\gamma^n}&=& \sum_{1}^{\Lambda} d_i(\gamma^n v_i(\e),v_i(\e))
\end{eqnarray*}
shows that for at least one  $i$ 
the sequence $d_i(\gamma^n v_i(\e),v_i(\e))$
must be unbounded.
Since $\gamma$ acts as an isometry, for any other vertex $v\in \mathcal{T}_i$: 
\begin{eqnarray*}
d_i(\gamma^n v_i(\e),v_i(\e))&\leq&d_i(\gamma^n v,v)+
d_i(\gamma^n v,\gamma^n v_i(\e))+d_i(v,v_i(\e))\\
&\leq&d_i(\gamma^n v,v)+2d_i(v,v_i(\e)).
\end{eqnarray*}
We infer that $\gamma$ does not stabilise  any 
finite set of vertices of $\mathcal{T}_i$.
In particular, $\gamma$ acts without inversion.
Now Proposition 25 of \cite{Se} implies the assertion.
\end{proof}
\begin{rem}
If  $m(s,t) < \infty$ for all $s,t \in S$, then the Coxeter group
itself has property FA of Serre.
(Concerning property FA, see {\rm\cite{Se}}, ex. 3 ,  p. 66.)
\end{rem}
Denote $I_1$, $I_2$ and $I_3$ the set of indices $i\in \{1,\ldots,\Lambda\}$
such that the corresponding trees have only one edge, only vertices of
valencies at most two, at least one vertex of valency at least three, 
respectively.
Lemma \ref{lem:trans} shows that the existence of one vertex of valency 
two implies that the tree contains an infinite axis of a translation
of amplitude two.
\par
Clearly, if $G$ is finite then $I_2=I_3=\emptyset$. On the other hand, if 
we denote $H_i \lhd \Gamma$ the intersection of the kernels
of the homomorphisms $\pi_j : \gamma \rightarrow \Aut(\mathcal{T}_j) \, , j\notin I_i$,
then, whenever $G$, equivalently $\Gamma$,
is infinite we have that $H_1=\{\e\}$ is trivial, $H_2$ is a solvable, normal subgroup 
of $\Gamma$ and,
if not trivial, $H_3$ contains non-abelian free subgroups.
\par
That $H_2$ is solvable follows from the facts that the
group of automorphisms of a tree of degree two is just
$\ZZ_2 \ltimes\ZZ$ and that the set $\pi_j,\ j\in I_2$
separates the points of $H_2$. indeed, $H_2$ embeds as a subgroup in a
direct sum of solvable groups. 
\begin{prop}%
\label{prop:trans}%
Let $\mathcal{T}$ be a tree with at least one vertex of valency at least three.
Assume that every pair of adjacent edges $e_1 = \{y,x\},\ e_2 = \{x,z\}$
defines a translation $u=u(e_1,e_2)$ on $\mathcal{T}$ with $uy=z$.
If $h_1, \ldots, h_l$  are non-trivial 
translations 
of $\mathcal{T}$,
then there is 
a pair of adjacent edges, defining a translation $v$,
such that for each $j \in \{1,\ldots, l\}$
the group generated by 
$h_j$ and $v$ in $\Aut(\mathcal{T})$ 
is isomorphic to the
free product $\ZZ \ast  \ZZ$.
\end{prop}
\begin{proof}
We shall use Klein's table tennis criterion,
of which a well fitting formulation, for our needs,
 can be found in \cite{BeCoHa95}
(Lemma 4.1).
\par
Since a group is acting on the tree there can at most be two valencies
of vertices. More precisely each vertex either is of valency at least three
or has a neighbour of this kind.
It is known that the boundary of the tree is infinite.
\par
By assumption, each $h_j$ has an attracting boundary point $b_{j}^+$
and a repulsing one $b_{j}^-$, which are connected by the axis
$a_j$ of the respective translation. We take a boundary point 
not contained in $\{b_{1}^+,\ldots,b_{l}^+\} \cup \{b_{1}^-,\ldots,b_{l}^-\}$
and find on a straight path towards this point a
vertex $x$ of valency at least three not belonging to
one of the axes $a_1,\ldots ,a_l$. 
The translations $\{h_1, \ldots, h_l\}$ do not fix this vertex and
our proof  is finished by the following well-known argument: 
\par 
Let $e$ be an  edge adjacent to $x$, 
but not belonging to one of the geodesics from $x$ to $h_ix$ neither to those
from $h^{-1}_ix$ to $x$, $i = 1,\ldots,l$.
We split the tree in two disjoint trees cutting  this edge. Let $V$
denote that part not containing the above geodesics and $U$ the one
containing them.
By assumption we find a translation $u$ moving $x$ into $V$.
The vertex $ux$ is adjacent to two different edges, which lie inside
$V$, since as an 
image of $x$ it has valency greater two.
Again by our assumption we find a translation $v$
whose axis entirely lies in $V$ and contains $ux$.
\par 
Now, for $h\in\{h_1,\ldots,,h_l\}$  and $j\in \ZZ \setminus \{0\}$
it is clear that  $h^j V \subset U$ and
on the other hand $v^j U \subset V$.
Since $v$ is of infinite order, Klein's criterion
implies that the group $<h,v>$ generated by $v$ and $h$ in
the group of automorphisms of $\mathcal{T}$
is the free product
$<h> \ast\ \ZZ = \ZZ \ast \ZZ$.
\end{proof}
\section{Factoriality}
We consider again the geometrical representation
$\sigma : G \rightarrow Gl(E)$. Associated to $G,S$ there is
the bilinear form $B:E\times E \rightarrow \RR$ whose matrix, with 
respect to the standard unit vectors of $E = \oplus_{s\in S} \RR e_s$,
has entries:
\[
B(e_s,e_t) = \left\{ \begin{array}{lcl}
1& \mbox{ if }& s=t\\
 -\cos(\pi \slash m(s,t))& \mbox{ if }&m(s,t)<\infty\\
 -1 & \mbox{ if }&m(s,t)=\infty . 
\end{array}\right.
\]
We shall call $(G,S)$ decomposable, if there exist non-empty subsets
$S_1,S_2 \subset S$,  such that $s,t\in S $ commute whenever
$s\in S_1$ and $t\in S_2$ or equivalently $B(e_s,e_t)=0$
and indecomposable otherwise.
\par
It is well known that the (left) regular representation of
a discrete group is factorial exactly if the group is icc, 
that is, the conjugation class of any group element different from the identity
is infinite. Our proof of this for a certain class of Coxeter groups relies 
very much on an irreducibility lemma of de~la~Harpe
for finite index subgroups of Coxeter groups (\cite{Harpe_87}, Lemma 1). 
It is not immediately clear that the complexified 
representations remain irreducible and we shall provide a proof.
\begin{prop}%
\label{prop:icc}
Let $(G,S)$ be an indecomposable Coxeter system, with $G$ infinite.
If the associated bilinear form $B$ is indefinite and non-degenerate
then $G$ is an icc-group.
\end{prop}
\begin{proof}
For $w\in G$ denote $C(w)$ its centraliser. The conjugation class
of $w$ is finite if and only if the index of $C(w)$ in $G$ is finite.
By Lemma 1 of \cite{Harpe_87} the image $\sigma(C(w))$ acts irreducible
on $E$ in this case. 
By Schur's lemma the commutant $\sigma(C(w))'$
is a % 
division algebra
over $\RR$. As it is finite dimensional it is isomorphic to 
$\RR$, $\CC$, or $\HH$.
We claim that for any $u\in G$, with $\sigma(u)\in \sigma(C(w))'$
the operator $\sigma(u)$ is a real multiple of the identity.
\par
The claim implies the proposition, because then it follows
that $\sigma(w)$, obviously an element of $\sigma(C(w))'$, 
commutes with all elements from $\sigma(G)$.
Since $\sigma$ is a faithful representation of $G$ we conclude
that $w$ is in the centre of the Coxeter group, but 
this is trivial (see \cite{Hum_90}, section 6.3). 
\par 
To establish the claim it suffices to show that any 
such $\sigma(u)$ has real spectrum.
If $\xi + \ii \eta \in \spec[u],\ \mbox{with } \xi,\ \eta \in \RR$, then by \cite{BonsallDuncan}, chap.1 Theorem 8,
$(\xi - \sigma(u))^2 + \eta^2$ is singular and hence equals $0$.
\par
If $\xi = 0$, then $\sigma(u^2)=\sigma(u)^2=- \eta^2$.
Since $\det{\sigma(u^2)} \in \{+1,-1\}$ we would have $\eta^2=1$
and $\sigma(u^2)=-1$, which is impossible in an infinite
Coxeter group.
\par
Now from $(\xi - \sigma(u))^2 + \eta^2 = 0$ we see
$2\xi  = \sigma(u^{-1})(\sigma(u^2) + \xi^2 +\eta^2)$.
Here the adjoint to the right hand side leaves
invariant the Tits cone $U$.
Hence $\xi$ must be strictly positive.
\par
We conclude that any $u$ with $\sigma(u)\in \sigma(C(w))'$
has its spectrum in the right half plane $\{z \ : \Re{z} >0\}$.
If some $z\in \spec[u]$ had  non-vanishing imaginary part then,
on the one hand, $z^k$, for some $k\in \NN$,   has negative real part, but 
on the other hand, by 
the spectral mapping theorem, $z^k$ is an element of the spectrum of
$\sigma(u)^k=\sigma(u^k)\in \sigma(C(w))'$.
\end{proof}
\begin{cor}
A Coxeter group as in the above proposition is not a
finite extension of an abelian group.
\end{cor}
\par
Let $\sigma_{\CC}: G \rightarrow \Gl(E\otimes_{\RR} \CC)$ be
the complexification of the geometric representation,
i.e. $\sigma_{\CC}(g) = \sigma(g)\otimes_{\RR} \id_{\CC}$, for all $g\in G$,
and extend $B$ canonically to a bilinear form $B_{\CC}$, which
clearly remains non-degenerate, if $B$ is.
\begin{lem}%
\label{lem:irred}
Suppose that  $G$ is infinite and $B$ non-degenerate.
Every subgroup of finite index in $G$ acts, by $\sigma_{\CC}$, irreducibly
on $E\otimes_{\RR} \CC$.
\end{lem}
\begin{proof}
We shall follow the first part of the arguments of de~la~Harpe.
We may assume that
$\Gamma$ is normal in $G$. Assuming $L_1 \neq E\otimes_{\RR} \CC$ to be a
non-trivial $\sigma_{\CC}(\Gamma)$-invariant subspace, one finds a generator
$s\in S$ such that $L_1\cap \sigma_{\CC}(s) L_1 = \{ 0 \}$. 
\par
The complex codimension one subspace
$H_s=\oplus_{s'\in S\setminus\{s\}}\CC e_{s'}$
is stabilised by $\sigma_{\CC}(s)$ and it has non-trivial
intersection with $L_1\oplus\sigma_{\CC}(s) L_1$, because the dimension of the
latter subspace is at least two. On the other hand 
 $L_1$ intersects trivially with $H_s$ since 
$\sigma_{\CC}(s)$ does not fix any of its non-zero elements.
We conclude that $L_1$ complements $H_s$ and is one-dimensional.
Especially we can write : $e_s = v + h$ for some $v\in L_1$
and some $h \in H_s$.
Now 
\begin{eqnarray*}
-e_s&=& \sigma_{\CC}(s)e_s \, = \, \sigma_{\CC}(s)v + h
\end{eqnarray*}
Subtracting, we see that $e_s= \frac{1}{2}(v-\sigma_{\CC}(s)v)\in
L_1\oplus\sigma_{\CC}(s) L_1$.
\par
The $G$-orbit $\mathcal{L}=\{L_1, \ldots , L_N\}$ of $L_1$
is finite since $\Gamma$ is of finite index in $G$,
and all of those complex lines are $\Gamma$-invariant,
by normality of $\Gamma$.
As $\dim L_j = 1$ there 
exist  homomorphisms $\lambda_j : \Gamma \rightarrow \CC^*$
by which $\Gamma$ acts.
Because $G$ acts irreducibly,
(notice that the extension to complex scalars is included 
in the Corollaire of chap.V \S 4 sec. 7 of \cite{Bourbaki_68})
on $E\otimes_{\RR} \CC$
the $G$-invariant sum $\oplus_{j=1}^N L_j$ equals the whole space.
Since $\sigma_{\CC}$ is a faithful representation
we conclude that $\Gamma$ is abelian, a contradiction
to the above corollary.
\end{proof}
\begin{rem}%
\label{rem:centr}
As in the proof of Proposition~\ref{prop:icc} one sees
that under the conditions of the above lemma a finite index
subgroup of $G$ has a trivial centraliser.
\end{rem}
\begin{prop}%
\label{prop:solv}
If $G$ is infinite and $B$ non-degenerate
then every  torsion free normal subgroup $\Gamma$ of finite index in $G$
contains no
non-trivial solvable normal subgroup.
\end{prop}
\begin{proof}
Let $H\lhd \Gamma$ be a solvable, normal subgroup of $\Gamma$ 
as in the statement of the 
proposition. 
We denote $\overline{H}^Z$ and $\overline{\Gamma}^Z$
the Zariski closures of $\sigma_{\CC}(H)$ respectively
of $\sigma_{\CC}(\Gamma)$ in $\Gl(E\otimes_{\RR} \CC)$. 
Clearly, $\overline{H}^Z$ is
a normal divisor of $\overline{\Gamma}^Z$ and, moreover
the connected component $\overline{H}^0$ (in the Zariski topology) of the 
identity in $\overline{H}^Z$ on the one hand is still
normal in $\overline{\Gamma}^Z$ and on the other hand
is of finite index in $\overline{H}^Z$. We claim 
that it reduces to the identity. This claim proves the proposition, 
since it implies that $\overline{H}^Z$ and hence $H$ are finite groups.
Because $\Gamma$ is torsion-free this is possible
only if $H=\{ \e \}$.
\par
The solvable Zariski connected group $\overline{H}^0$
has a common eigenvector $v\in E\otimes_{\RR} \CC$, 
as follows from the Lie-Kolchin Theorem, see e.g.\
Corollary 10.5 in \cite{Borel_91}.
Therefore, there exists a character (a continuous multiplicative function)
$\alpha_v :\overline{H}^0  \rightarrow \CC^{\ast}$
from $\overline{H}^0$ to the multiplicative group of $\CC$,
such that $h\ v = \alpha_v(h)\ v$.
Since $\overline{H}^0$ is normal in $\overline{\Gamma}^Z$ any
vector $\gamma v$, with $\gamma \in \overline{\Gamma}^Z$,
is also a common eigenvector, the corresponding
character is $\alpha_{\gamma v}(.)=\alpha_v(\gamma^{-1}.\ \gamma)$.
\par
Let 
$\mathcal{V}=\{u\in E\otimes_{\RR} \CC : h\ u = \alpha_u(h) u \mbox{ for some }\alpha_u
\mbox{ as  above }\}$. This set is $\overline{\Gamma}^Z$ invariant
and spans an $\overline{\Gamma}^Z$-invariant subspace. 
We have seen that it is non-trivial, and by the irreducibility
of $\sigma_{\CC}(\Gamma)$ it must equal $E\otimes_{\RR} \CC$.
\par
Now the trace of $B_{\CC}$ is positive and $\mathcal{V}$
contains a basis, hence for some $u\in\mathcal{V}$
$B_{\CC}(u,u)\neq 0$. From
\begin{eqnarray*}
\alpha_u(h)B_{\CC}(u,u) &= & B_{\CC}(hu,u) \, = \, B_{\CC}(u,h^{-1}u)\\
&=&\alpha_w(h^{-1})B_{\CC}(u,u)\quad \forall h\in\overline{H}^0
\end{eqnarray*}
we infer $\alpha_{u}\in \{+1, -1\}$.
\par
As above we conclude that the set
\[\mathcal{V}_1=\{u\in E\otimes_{\RR} \CC : h\ u = \alpha_u(h) u, \;%
\alpha_u:\overline{H}^0\to \{+1, -1\}
\mbox{ is a character }\}\]
contains a basis. With respect to one such basis the elements of 
$\overline{H}^0$ consist of diagonal matrices with entries from $\{+1, -1\}$.
Since $\overline{H}^0$ is Zariski connected it must be trivial.
\end{proof}
\section{Simplicity of the Regular $C^{\ast}$-Algebra.}
\begin{thm}
If $(G,S)$ is an indecomposable Coxeter system, with $G$ infinite, such that
the associated bilinear form $B$ is indefinite and non-degenerate
then its (left) regular $C^{\ast}$-algebra is simple with unique trace.
\end{thm}
Before proving the theorem we shall
establish a lemma:
\begin{lem}
If $(G,S)$ is as in the theorem then all trees $\mathcal{T}_1,\ldots,
\mathcal{T}_{\Lambda}$ have vertices of valency at least three.
\end{lem}
\begin{proof}
Assume that $\mathcal{T}_i$ has only vertices of valency two.
From Lemma~\ref{lem:morph} we see that for $j$ in the orbit
$\mathcal{O}=\{\pi(\dot{g})i : \dot{g}\in G/\Gamma\}$
the tree 
$\mathcal{T}_j$ is isomorphic to $\mathcal{T}_i$, hence has only vertices
of valency two.
\par
The group $W$ generated by the reflections $\cup_{j\in\mathcal{O}} T_j$
contains $\Gamma':= W\cap \Gamma$ as a
normal torsion-free subgroup of finite index, which is normal in $G$ too.
The set of homomorphisms
$\pi_j | \Gamma' : \Gamma' \rightarrow \Aut (\mathcal{T}_j)$, $ j\in
\mathcal{O}$
is faithful, since the assumption that $\pi_j(\gamma)=\id$ 
for all $ j\in \mathcal{O}$
implies that $\gamma$ acts, by conjugation, trivially on all
reflections in $\cup_{j\in\mathcal{O}} T_j$, i.e.\ on all
reflections in $W$. From this and from the fact that $\gamma \in W$ we have
that $\gamma$ is in the centre of $W$. But the centre is trivial, 
since $W$ is infinite.
\par
We infer that $\Gamma'$ is solvable. By Proposition \ref{prop:solv}
it is trivial and $W$ finite. This is a contradiction. 
\end{proof}
\par\noindent
\begin{proof}[Proof of the theorem:]
Let $\Gamma$ be a torsion free, normal subgroup of finite index in $G$.
Since $G$ is an icc-group the results of Bekka and de~la~Harpe 
\cite{BekkaHarpe_2000}
show that it suffices to prove the assertions for $\Gamma$.
We shall use the concept of weak Powers groups in
the sense of Boca and Nitica \cite{BocaNitica_88}.
\par
Let $F\subset C_h$ be a finite subset of the $\Gamma$-conjugation-class
of some $h\in\Gamma$.
Lemma \ref{lem:translation} shows that there is a tree $\mathcal{T}$ on which
$h$, and hence all elements of its conjugation-class act as translations.
This tree has a vertex of valency at least three.
We may apply Proposition \ref{prop:trans} and find $v\in \Gamma$ 
such that for any $k\in F$ the 
subgroup $<k,v>$ generated by them is isomorphic to $\ZZ \ast \ZZ$.
\par
The proof of Lemma~2.2 of \cite{BeCoHa95} shows that 
we find a constant $C>0$ and 
$v\in \Gamma$ such that
for all $k\in F$:
\begin{eqnarray}%
\label{eq:2estimate}
\norm{\sum_{j=1}^{\infty} a_j \lambda_{\Gamma}(v^{-j}k v^{j})}&\leq& C \ %
\norm[2]{a} \qquad \forall a \in l^2(\ZZ^+).
\end{eqnarray}
\par
Armed with this, the computations in the proof of 
Lemma 2.2 in \cite{BocaNitica_88} prove the following fact which we state as a lemma.
\begin{lem}%
\label{lem:estimate}
Given a finite linear combination
$x=\sum_{k\in F} a_k k \in \CC\Gamma$
with $\e\notin F$ and $\epsilon>0$
there exist $n\in\NN$ and $v_1,\ldots, v_n \in \Gamma$ 
such that in $C^{\ast}_{\lambda}(\Gamma)$:
\begin{eqnarray}%
\label{eq:estimate}
\norm{\frac{1}{n}\sum_{j=1}^n%
  \lambda_{\Gamma}(v_j)\lambda_{\Gamma}(x)\lambda_{\Gamma}(v_j)^{\ast} }
&\leq & \epsilon .
\end{eqnarray}
\end{lem}
The arguments in the proof of Lemma~2.1 in \cite{BeCoHa95} show that
$C^{\ast}_{\lambda}(\Gamma)$ is simple.
Finally, the
uniqueness of the trace is an immediate consequence of 
inequality (\ref{eq:estimate}).
\end{proof}
\begin{rem}
It is not hard to see, that $\Gamma$ is a weak Powers group.
By Remark~\ref{rem:centr}, the centraliser of $\Gamma$ in $G$ is trivial, 
hence the Coxeter group itself is an ultra-weak Powers group in the sense of 
B\'{e}dos~\cite{Bedos_91}.
\end{rem}
\bibliographystyle{amsplain}
\providecommand{\bysame}{\leavevmode\hbox to3em{\hrulefill}\thinspace}

\vspace{\fill}
{{\bf Author's address}\\
Gero Fendler\\
Finstertal 16\\
D-69514 Laudenbach\\
Germany\\
\small
e-mail: gero.fendler@t-online.de}
\end{document}